\documentclass[pdflatex,sn-mathphys-num]{sn-jnl}% Math and Physical Sciences Numbered Reference Style
\usepackage{graphicx}%
\usepackage{multirow}%
\usepackage{amsmath,amssymb,amsfonts}%
\usepackage{amsthm}%
\usepackage{mathrsfs}%
\usepackage[title]{appendix}%
\usepackage{xcolor}%
\usepackage{textcomp}%
\usepackage{manyfoot}%
\usepackage{booktabs}%
\usepackage{algorithm}%
\usepackage{algorithmicx}%
\usepackage{algpseudocode}%
\usepackage{listings}%
\theoremstyle{thmstyleone}%
\newtheorem{theorem}{Theorem}%  meant for continuous numbers
\newtheorem{proposition}[theorem]{Proposition}%

\newtheorem{corollary}[theorem]{Corollary}

\theoremstyle{thmstyletwo}%

\theoremstyle{thmstylethree}%
\newtheorem{definition}{Definition}%

\begin{document}

\title[Characterizations of Faces]{Characterizations of Faces of Convex Sets in Infinite-dimensional Vector Spaces}


\author{\fnm{Valentin V.} \sur{Gorokhovik}}\email{gorokh@im.bas-net.by}


\affil{\orgdiv{Institute of Mathematics}, \orgname{The National Academy of Sciences of Belarus}, \orgaddress{\street{Surganov Str., 11}, \city{Minsk}, \postcode{220012}, %\state{State},
\country{Belarus}}}

%\affil[2]{\orgdiv{Department}, \orgname{Organization}, \orgaddress{\street{Street}, \city{City}, \postcode{10587}, \state{State}, \country{Country}}}
%
%\affil[3]{\orgdiv{Department}, \orgname{Organization}, \orgaddress{\street{Street}, \city{City}, \postcode{610101}, \state{State}, \country{Country}}}


\abstract{In the paper three different characterizations of faces of convex sets, belonging to infinite-dimensional real vector spaces, are presented. The first one is formulated in the terms of generalized semispaces, the second --- in the terms of compatible complete (total) preorders, and the third --- in the terms of step-affine functions. All three characterization are equivalent each other and extend to infinite-dimensional vector spaces the lexicographical characterization of faces established in finite-dimensional settings by Martinez-Legaz J.-E. (Acta Mathematica Vietnamica. 1997. Vol. 22, No.~1, P. 207--211).
}




\maketitle

\section{Introduction}\label{sec1}

Throughout this paper $X$ is nontrivial ($X \ne \{0\}$) real vector space which is supposed to be infinite-dimensional, unless otherwise stated. Besides, it is not supposed that some other structures (in particular, the topological one) are defined on $X$.

A set
$Q \subseteq X$ is called \cite{Rock,Boris} {\it convex}, if for any $x,y \in Q$ the segment $[x,y]:=\{\alpha x + (1-\alpha)y \mid \alpha \in [0,1]\}$ is entirely contained in $Q$. %The empty set $\varnothing$ is considered convex by definition.

A subset $K \subseteq X$ is called a \textit{cone} if $\lambda x \in K$ for all $x \in K$ and all $\lambda >0$.

In other words a subset $K \subseteq X$ is a cone if along with each point $x \in K$ the subset $K$ contains the entire ray emanating from the origin and passing through the point $x$.
Note, that the origin not necessarily belongs to a cone $K$.

A cone $K \subseteq X$ is convex if and only if $x + y \in K$ whenever $x,y \in K$ or, equivalently, if and only if $K + K \subseteq K$.

A cone $K \subseteq X$ is said to be \textit{asymmetric}, if $K \bigcap (-K) = \varnothing$. A convex cone $K \subset X$ is asymmetric if and only if $0 \not\in K$.
%Рассмотрим специальный подкласс выпуклых множеств, называемых полупространствами.

%Выпуклое подмножество $H \subset X$ векторного пространства $X$ называется \cite{Las} \textit{выпуклым полупространством в $X$} или \cite{GS98,GS2000,Gor20} просто  \textit{полупространством в $X$}, если его дополнение $X \setminus H$ также является выпуклым множеством. Разумеется, что если $H$ полупространство в $X$, то и $X \setminus H$ также является полупространством в $X$.
%
%Полупространство, которое является конусом, называется \textit{коническим полупространством}.

A nonempty convex subset $F$ of a convex set $Q \subset X$ is called (see, for instance, \cite{Rock,Boris}) \textit{a face} of the set $Q$ if it satisfies the following property: if for some $u,v \in Q$ there exists $\alpha \in (0,1)$ such that the point $\alpha u + (1-\alpha)v$ is contained in $F$ then $u,v \in F$.

In other words, a nonempty convex subset $F \subset Q$ is \textit{a face} of $Q$ if every segment of $Q$, having in its relative interior an element of $F$, is entirely contained in $F$. The set $Q$ itself is its own face and the empty set is considered as a face of any convex set $Q$. A nonempty face $F$ of a set $Q$ is called \textit{proper} if it is not equal to the set $Q$ itself.

It is easy to see that when $K$ is a convex cone, each face $F$ of $K$ is a convex cone too. Moreover, when $K$ is an asymmetric convex cone, each face $F$ of $K$ is an asymmetric convex cone as well. An intersection of any family of faces of a convex set $Q$ is also a face of $Q$ while the union of a family of faces of a convex set $Q$ is not always a face of $Q$. However, if the family of faces is linearly ordered by the inclusion relation, then the union of such a family is also a face of $Q$. If $F$ is a face of $Q$ while $G$ is a face of $F$, then $G$ is also a face of $Q$. These statements are fairly obvious and can be found in most works dealing with faces (see, for example, \cite{Rock,Boris,Bronsted,Millan}).

A subset $H_=(l,\alpha):=\{x \in X \mid l(x) = \alpha\}$, where $l: X \to {\mathbb R}$ is a nonzero linear function defined on $X$, and $\alpha$ is a real number. is called \textit{a hyperplane} in $X$. A hyperplane $H_=(l,\alpha)$ is said to be \textit{support} to a convex set $Q$, if $Q \subset H_\geq (l,\alpha):=\{x \in X \mid l(x) \geq \alpha\}$ and $H_=(l,\alpha) \bigcap Q \ne \varnothing$. If a hyperplane $H_=(l,\alpha)$ is support to a convex set $Q$, then the intersection $H_=(l,\alpha)$ with the set $Q$, i.e. the subset $H_=(l,\alpha) \bigcap Q$, is a face of $Q$. Such faces of a convex set $Q$ is called  \cite{Rock} \textit{exposed}.
In other words, a face $F$ of a convex set $Q \subset X$ is exposed, if it coincides with the set of points from $Q$, at which  a linear function $l: X \to {\mathbb R}$ attains its minimum over the set $Q$, i.e. if $F=\{x \in Q \mid l(x) = {\rm min}_{y \in Q}l(y)\}$, where $l: X \to {\mathbb R}$ is a nonzero linear function. Not every face is exposed (an example of a non-exposed face can be found in \cite[p.~179]{Rock}). If $F$ is an exposed face of $Q$, and $G$ is an exposed face of $F$, then $G$ is a face of $Q$, but not necessarily an exposed face of $Q$.

\begin{definition}\label{def1} A proper face $F$ of a convex set $Q \subset X$ will be called \textit{lexicographically exposed}, if there exists a finite ordered family $\{l_1,l_2,\ldots,l_m\}$ of nonzero linear functions defined on $X$ such that
\begin{equation}\label{e0.1}
F = \{x \in Q \mid l_i(x) = {\rm min}_{y \in F_{i-1}}l_i(y), i = 1,2,\ldots,m\},
\end{equation}
where $F_0 = Q, F_i = \{x \in F_{i-1} \mid l_i(x) = {\rm min}_{y \in F_{i-1}}l_i(y),\},  i = 1,2,\ldots,m$.
\end{definition}

Without loss of generality, we can assume that the family $\{l_1,l_2,\ldots,l_m\}$ of linear functions in the definition of a lexicographically exposed face is linear independent.

For exposed faces the equality \eqref{e0.1} holds under $m = 1$ and, consequently, each exposed face is lexicographically exposed. Note also, that each subset $F_i, i = 1,2,\ldots,m,$ in \eqref{e0.1} is an exposed face of the preceding face $F_{i-1}$.

Informally, a lexicographically exposed face $F$ of a convex set $Q$ can be characterized as a face that can be get in the result of the following finite sequential procedure. Choose an exposed face $F_1$ of the set $Q$ that contains $F$, if $F = F_1$ the procedure is complete. In the case of $F \ne F_1$ we choose an exposed face $F_2$ of the face $F_1$ that contains $F$, if $F = F_2$ then the procedure is complete. In the case of $F \ne F_2$ we continue the procedure further. If after a finite number of steps, say $m$ steps, we get an exposed face $F_m$ of the previous face $F_{m-1}$ such that $F_m = F$, then the face $F$ is lexicographically exposed.

As it was proved by J.-E. Martinez-Legaz in \cite[Theorem 2]{ML} in finite-dimensional vector spaces each proper face of any convex set is lexicographically exposed.
The authors of the paper \cite{Millan} put the following question: is it possible to reasonably generalize the lexicographical characterization of faces\footnote{In \cite{Millan} lexicographically exposed faces are called \textit{eventually exposed}.} to convex sets belonging to infinite-dimensional vector spaces?

In the present paper we present three different (but, of course, equivalent) characterizations of arbitrary faces of convex sets in infinite-dimensional spaces (see Theorems \ref{thA}, \ref{thB}, \ref{thC}), which in the finite-dimensional setting are equivalent to the lexicographical characterization of faces established by J.-E.~Martinez-Legaz in \cite[Theorem~2]{ML}.

%%%%%%%%%%%%%%%%%%%%%%%%%%%%%%%%%%%%%%%%%%%%%%%%%%%%%%%%%%%%%%%%%%%%%%%%%%%%%%%%%%%%%%%%%%%%%%%%%%%%%%%%%%%%%%%%%%%%%%%%%%%

\section{Some properties of faces of convex sets.}

Let's begin with a proposition in which we will present some properties of the faces of convex sets that we will need in what follows.

\begin{proposition}\label{pr1}
Let $Q$ be a convex set in a real vector space $X$ and let $F$ be a proper face of the set $Q$.
Then

$($i$)$ the set $Q \setminus F$ is convex; %and, in addition, ${\rm icr}(Q \setminus F) = {\rm icr}Q$;

$($ii$)$ for arbitrary points $x \in F$ and $y \in Q \setminus F$ the open halfline $\{x - \tau(y-x) \mid \tau >0\}$ does not intersect the set $Q$, i.e. $\{x - \tau(y-x) \mid \tau >0\} \cap Q = \varnothing$;

$($iii$)$ ${\rm aff}F \bigcap (Q\setminus F) = \varnothing$ and hence ${\rm aff}F \bigcap Q = F$.
\end{proposition}

\textit{Proof.}
(\textit{i}) Let $x,y \in Q \setminus F$. Assume that there exists a point $z \in F$ such that $z \in (x,y)$. But then by the definition of a face we have $x,y \in F$. However, this contradicts the choice of $x$ and $y$, and hence $[x,y] \subset Q \setminus F$.

($ii$) Arguing by contradiction, suppose that for some $t > 0$ the point $x - t(y-x)=: z$ belongs to $Q$. Then $\left(1-\displaystyle\frac{1}{1+t}\right)y + \displaystyle\frac{1}{1+t}z = \left(1-\displaystyle\frac{1}{1+t}\right)y + \displaystyle\frac{1}{1+t}(x -t(y-x)) = x \in (y,z)$, and $y,z \in Q$. Since $F$ is a face of $Q$, then $y,z \in F$, but this contradicts the choice $y \in Q \setminus F$.

($iii$) The assertion is trivial when $F$ is singleton. Assume that $F$ contains more than one point. Since $F$ is convex, then ${\rm aff}F = \{\tau x_1 + (1-\tau)x_2 \mid x_1,x_2 \in F\,\,\text{and}\,\,\tau  \in {\mathbb R}\}$. Suppose on the contrary to the assertion ($iii$) that ${\rm aff}F \bigcap (Q\setminus F) \ne \varnothing$ and take a point $y \in {\rm aff}F \bigcap (Q\setminus F)$.  Since $y \not\in F$, then $y = t x_1 + (1-t)x_2$ for some $x_1,x_2 \in F, x_1 \ne x_2,$ and some $t>1$. The points $y,x_1,x_2$ lies on the same line $l(x_1,x_2):= \{\tau x_1 + (1-\tau)x_2 \mid \tau  \in {\mathbb R}\}$, and $x_1 \in (y, x_2)$. It implies that $x_2$ belongs to the open halfline $\{x_1 - \tau(y-x_1) \mid \tau >0\}$. But this contradicts the assertion ($ii$) since $x_2 \in F \subset Q$, \hfill $\square$

Note that the property ($i$) of Proposition \ref{pr1} was also noted in \cite[Proposition~4.7]{Millan}.

\smallskip

\section{The characterization  of faces of convex sets by halfspaces.}

A convex set $H \subset X$ of a vector space $X$ is called \cite{Las} \textit{a convex halfspace in $X$} or \cite{GS98,GS2000,Gor20} simply  \textit{a halfspace in $X$}, if its compliment $X \setminus H$ is also a convex set. A halfspace, which is a cone, is called \textit{a conical halfspace}. Clearly, that a complement of any (conical) halfspace is also a (conical) halfspace. Halfspaces which are distinct from the empty set $\varnothing$ and the space $X$  are called \textit{proper}.

A convex set $H \subset X$ is a halfspace if and only if \cite[Theorem~3]{GS98} its recession cone $0^+H := \{y \in X \mid x +
ty \in H\,\,\text{for all}\,\,x \in H\,\,\text{and all}\,\,t>0\}$ is a conical halfspace.

An asymmetric convex cone $K \subset X$ is a conical halfspace if and only if \cite[Corollary, p.~165]{GS98} the set $L := X \setminus (K\bigcup(-K))$ is a vector subspace in $X$.

A pointed ($K \bigcap (-K) \ne \varnothing$) convex cone $K \subset X$ is a conical halfspace if and only if \cite[Theorem~2]{GS98} $X = K \bigcup (-K)$, and its asymmetric part  $K^>:=K \setminus (-K)$ is also a halfspace in $X$, because $X \setminus (K^> \bigcup (-K^>)) = K \bigcap (-K)$ is a vector subspace in $X$.

In infinite-dimensional vector spaces there are three, and in finite-dimensional spaces two, different types of halfspaces (for details on the types of halfspaces, see \cite{GS98,GS2000,Gor20}), according to which all halfspaces are divided into three (in finite-dimensional spaces into two) disjoint classes. One of these classes consists of the so-called generalized semispaces generated by various affine manifolds from $X$, another class consists of the complements of generalized semispaces, and the halfspaces of the third class, which exist only in infinite-dimensional spaces, are called gap-halfspaces\footnote{In \cite{GS98,GS2000,Gor20} the following terminology is used to name the types of halfspaces: generalized semispaces are called non-pointed Dedikin halfspaces, their complements are called pointed Dedikin halfspaces, and gap halfspaces are called non-Dedekind halfspaces.}. An example of a gap-halfspace is given in \cite{GS98}.

Generalized semispaces are defined as follows.

Let $M \subset X$ be an affine manifolds from $X$. \textit{A generalized semispace}, generated by $M$, is a maximal (by inclusion) convex set, belonging to $X \setminus M$. The existence of generalized semispaces, generated by an arbitrary affine manifolds  follows from the Kuratowski-Zorn lemma \cite{Hille,Jech}. Moreover, it follows also from the Kuratowski-Zorn lemma, that for any convex set $C \subset X$ such that $C \bigcap M = \varnothing$, where $M \subset X$ is an affine manifolds from $X$, there exists an generalized semispace, generated by $M$, which contains $C$. Generalized semispaces are halfspaces in $X$. To show this we will need the following auxiliary assertion, characterizing the geometrical structure of generalized semispaces.

\begin{proposition}\label{pr2}
Let $S$ be a generalized semispace generated by an affine manifolds $M \subset X$ and let a point $a$ belong to $M$. Then $S = S_a :=  \{a + t(x-a) \mid x \in S, t >0\}$.
\end{proposition}

\textit{Proof.}
Let us show that for any $a \in M$ the set ${S}_a := \{a + t(x-a) \mid x \in S, t >0\}$ is convex, contains $S$, and does not intersect $M$.
The convexity of $S_a$ follows from the convexity of $S$. In addition, $S_a \bigcap M = \varnothing$. Indeed, otherwise there would exist a point $\bar{x} \in S$ such that the halfline $\{a + t(\bar{x}-a) \mid t >0\}$ would belong to $M$ and, consequently, (for $t=1$) the point $\bar{x}$ would belong to $M$, which contradicts $\bar{x} \in S$. Since $S$ is a maximal (by inclusion) convex set such that $S \bigcap M = \varnothing$, then $S_a = S$. \hfill $\square$

It is easy to see  from the definition of the set $S_a$, given in Proposition \ref{pr2}, that $\{a\} \bigcup S_a$ is also a convex set, and this in turn implies that $M \bigsqcup S$ (here and below $\bigsqcup$ means disjoint union) is a convex set.

If $S \subset X$ is a generalized semispace generated by an affine manifolds $M \subset X$, then the $S^- := \{x \in X \mid x = 2a - y, y \in S\}$, where $a$ is an arbitrary point from $M$, is also a generalized semispace generated by $M$, and $X = S^-\bigsqcup M \bigsqcup S$.

Indeed, assume that $X \ne S^-\bigsqcup M \bigsqcup S$ and let $\bar{x} \in X \setminus (S^-\bigsqcup M \bigsqcup S)$. Then for any $a \in M$ the halfline  $l^+(a,\bar{x}):=  \{a + t(\bar{x} - a) \mid t > 0\}$ belongs to $X \setminus (S^-\bigsqcup M \bigsqcup S)$ and hence  ${\rm conv}(S \bigcup l^+(a,\bar{x}))$ is a convex set which does not intersect $M$. In addition, $S$ is a proper subset of ${\rm conv}(S \bigcup l^+(a,\bar{x}))$. But this contradicts the maximality of $S$.

It follows from the equality $X = S^-\bigsqcup M \bigsqcup S$ and the convexity of the sets $S \bigsqcup M$ and $S^- \bigsqcup M$ that generalized semispaces $S$ and $S^-$ are halfspaces, and, consequently, the sets $H^+:= M \bigsqcup S$ and $H^-:= S^-\bigsqcup M$ (as complements of halfspaces) are also halfspaces in  $X$, however their type is distinct from the type of generalized semispaces.  %We will refer to the halfspaces $H^+:= M \bigsqcup S$ and $H^-:= S^-\bigsqcup M$ as \textit{pointed generalized semispaces}, generated by the affine manifold $M$.

%Будем говорить, что заостренное обобщенное семипространство $H^+:= M \bigsqcup S$ является \textit{опорным} к выпуклому множеству $Q \subset X$, если $Q \subset H^+$ и $Q \bigcap M \ne \varnothing$.

%Поскольку каждое обобщенное семипространство $S \subset X$ является выпуклым подмножеством (более того, полупространством) в $X$, то в соответствии с Определением 1 на $S$ может быть задано отношение доминирования $\unlhd_S$ и построено семейство открытых компонент ${\mathcal O}(S)$,  которое (см. Предложение \ref{pr3.7}) линейно упорядочено отношением $\unlhd_S^*$.
%
%Будем говорить, что обобщенное полупространство $S \subset X$ имеет \textit{конечный ранг}, если семейство его открытых граней ${\mathcal O}(S)$ конечно (определение ранга для произвольного полупространства дано в \cite{GS98,GS2000,Gor20}).

\begin{theorem}\label{thA}
Let $Q$ be a convex set in a real vector space $X$.
A subset $F \subset Q$ is a proper face of the set $Q$ if and only if in $X$ there exist a nonempty affine manifold $M \subset X$ and a generalized semispace $S$, generated by $M$, such that
$Q \subset M \bigcup S$ and $F = M \bigcap Q$.
\end{theorem}

\textit{Proof.}
\textit{Necessity}. Let $F \subset Q$ be a proper face of a convex set $Q \subset X$. It follows from the assertions ($i$) and ($iii$) of Proposition \ref{pr1} that the set  $Q \setminus F$ is convex and ${\rm aff}F \bigcap (Q\setminus F) = \varnothing$. Choose an arbitrary generalized semispace $S \subset X$, generated by the affine manifold ${\rm aff}F$, such that $Q \setminus F \subset S$. Then $Q \subset {\rm aff}F \bigcup S$ and $F = {\rm aff}F \bigcap Q$. Thus, the assertion of the necessary part of the theorem under consideration holds with $M = {\rm aff}F$ and any generalized semispace $S$, generated by ${\rm aff}F$, such that $Q \setminus F \subset S$.

\textit{Sufficiency}. Assume that $Q \subset M \bigcup S$ and $F = M \bigcap Q$, where $M$ is an affine manifold in $X$, and $S$ is a generalized semispace, generated  $M$. Since $S \bigcap M = \varnothing$ then $Q =(Q \bigcap S) \bigsqcup F$. Hence, for any $y,z \in Q$ the following three different cases are possible: 1) $y,z \in Q \bigcap S$, 2) $y \in Q \bigcap S, z \in F$, and 3) $y,z \in F$. Let us show that if the points $y,z \in Q$ are such that $(y,z) \bigcap F \ne \varnothing$, then the first two cases are impossible. Indeed, if $y,z \in Q \bigcap S$, then since the set $Q \bigcap S$ convex (as the intersection of two convex sets), then $[y,z] \subset Q \bigcap S$ and, consequently, $[y,z] \bigcap F = \varnothing$. Suppose now that $y \in Q \bigcap S, z \in F$ and $(y,z) \bigcap F \ne \varnothing$. Then the points $y,z$ and any point $u \in (y,z) \bigcap F$ lie on the same line $l(u,z) := \{\alpha u + (1 - \alpha)z \mid \alpha \in {\mathbb R}\}$, which belongs to $M$. This, however, contradicts the fact that $y \in Q \bigcap S$. Thus, if points $y,z \in Q$ are such that $(y,z) \bigcap F \ne \varnothing$, then $y,z \in F$ and, therefore, the subset $F = M \bigcap Q$ is a face of the set $Q$. \hfill $\square$
%%%%%%%%%%%%%%%%%%%%%%%%%%%%%%%%%%%%%%%%%%%%%%%%%%%%%%%%%%%%%%%%%%%%%%%%%%%%%%%%%%%%%%%%%%%%

\section{The characterization of faces of convex sets by compatible complete preorders.}

A preorder relation (a reflexive and transitive binary relation) $\preceq$ defined on $X$ is said \cite{Peressini} to be \textit{compatible} with algebraic operations of the space $X$ if
$$x \preceq y \Longrightarrow x+ z \preceq y+ z\,\,\forall\,\,x,y,z \in X$$
and
$$x \preceq y \Longrightarrow \lambda x \preceq \lambda y\,\,\forall\,\,x,y \in X, \lambda > 0.$$

If $\preceq$ is a compatible preorder relation then
$$x \preceq y \Longleftrightarrow y - x \in P_\preceq,$$
where $P_\preceq := \{x \in X \mid 0 \preceq x\}$ is  the (pointed, convex) cone of positive elements of the relation $\preceq$.

If $P_\preceq \bigcup (-P_\preceq) = X$, then the the preorder relation  $\preceq$ is total, i.e. such that for any $x,y \in X$ at least one of corelations $x \preceq y$ or $y \preceq x$ holds.

%Будем говорить, что \textit{ранг} согласованного отношения предпорядка $\preceq$ конечен, если векторное подпространство $P_\preceq \bigcap (-P_\preceq)$ имеет конечную коразмерность.
%Подмножество $H \subset X$ является полупространством в векторном пространстве $X$ тогда и только тогда, когда его рецессивный конус $0^+H := \{h \in X \mid x + th \in H\,\,\text{для всех}\,\,x \in H\,\,\text{и}\,\,t \ge 0\}$ является коническим полупространством. Так как $0^+H$ заостренный выпуклый конус и $0^+H \bigcup (-0^+H) = X$, то
%полупространство $H$ порождает на $X$ согласованное полное отношение предпорядка $\preceq_H$, определенное следующим образом:
%$$x \preceq_H y \Longleftrightarrow y - x \in 0^+H$$.

\begin{theorem}\label{thB}
Let $Q$ be a convex set in a real vector space $X$.
A subset $F \subset Q$ is a proper face of the set $Q$ if and only if
there exists a compatible total preorder relation $\preceq$ defined on $X$ such that $F = {\rm Min}(Q \mid \preceq) := \{ x \in Q \mid x \preceq y\,\,\text{for all}\,\,y \in Q\}$.
\end{theorem}

\textit{Proof.}
\textit{Necessity.}
It follows from the necessity part of Theorem \ref{thA} that for any proper face $F$ of a convex set $Q \subset X$ there exist an affine manifold $M \subset X$ and a generalized semispace  $S \subset X$, generated by $M$, such that $Q \subseteq M \bigcup S$ and $F = M \bigcap Q$.

%Как  уже было отмечено выше, выпуклое подмножество $H \subset X$ векторного пространства $X$ является полупространством в $X$ тогда и только тогда, когда его рецессивный конус $0^+H := \{y \in X \mid x + ty \in H\,\,\forall x \in H\,\,\text{и}\,\,\forall t > 0\}$ является коническим полупространством.

Since a generalized semispace $S$ is a halfspace then, as it was noted above, its recession cone $0^+S$ is a conical halfspace. Letting
$x \preceq_S y \Longleftrightarrow y - x \in 0^+S$ for any $x,y \in X$, we define on $X$ a compatible preorder relation $\preceq_S$,
which is total due to the equality $0^+S \bigcup (-0^+S) = X$.

Let us prove that $y \prec_S x$ ($y \prec_S x$ means that $y \preceq x$, but $x \not\preceq y$) for any $x \in S$ and any $y \in X \setminus S$. Arguing by contradiction, suppose that for some $\bar{x} \in S$ and $\bar{y} \in X \setminus S$, $\bar{y} \not\prec_S \bar{x}$. Since $\preceq_S$ is total, this assumption is equivalent to $\bar{x} \preceq_S \bar{y}$, which in turn is equivalent to the inclusion $\bar{y} - \bar{x} \in 0^+S$. By the definition of the recessive cone $0^+S$ we then have that $\bar{x} + (\bar{y} - \bar{x}) = \bar{y} \in S$, while by assumption $\bar{y} \in X \setminus S$. We have arrived at a contradiction, which proves the required statement.

Since $M \subset X\setminus S$, then $y \prec_S x$ for any $x \in S$ and any $y \in M$ and, therefore, $M \supset {\rm Min}(M \bigcup S \mid \preceq_S) := \{ x \in M \bigcup S \mid x \preceq y\,\,\text{for all}\,\,y \in Q\}$.

On the other hand, for any $a \in M$ the equality $M \bigcup S = a + 0^+S$ holds, from which it follows that $a \preceq_S x$ for any $x \in M \bigcup S$. Hence, $M \subset {\rm Min}(M \bigcup S \mid \preceq_S)$, and since the reverse inclusion also holds, $M = {\rm Min}(M \bigcup S \mid \preceq_S)$.

Taking into account  that $Q \subseteq M \bigsqcup S$ and $F = M \bigcap Q$, from the last equality we obtain $F = {\rm Min}(Q \mid \preceq_S)$.

\textit{Sufficiency.} Let $\preceq$ be a compatible total preorder relation on $X$ and let $F = {\rm Min}(Q \mid \preceq)$. Assume that the points $u,v \in Q$ are such that $(u,v) \bigcap {\rm Min}(Q \mid \preceq) \ne \varnothing$. Then there exists $\alpha \in (0,1)$ such that $\alpha u  + (1-\alpha) v \in {\rm Min}(Q \mid \preceq)$. Let us show that in this case $u$ and $v$ belong to  ${\rm Min}(Q \mid \preceq)$ as well. Indeed, if we would have $u \not\in {\rm Min}(Q \mid \preceq)$, then there would be a point $\bar{y} \in Q$ such that $\bar{y} \prec u$.
Since the relation $\preceq$ is compatible with the algebraic operations of the space $X$, it would follow from $\bar{y} \prec u$ that $\alpha \bar{y} + (1- \alpha)v \prec \alpha u + (1-\alpha) v$. Since $\alpha \bar{y} + (1- \alpha)v \in Q$, the last relation contradicts the fact that $\alpha u + (1-\alpha) v \in {\rm Min}(Q \mid \preceq)$. Consequently, $u,v \in {\rm Min}(Q \mid \preceq)$ and, hence, $F := {\rm Min}(Q \mid \preceq)$ is a face of the set $Q$.
\hfill $\square$

\smallskip

%Будем говорить, что определенное на векторном пространстве $X$ согласованное полное отношение предпорядка $\preceq$ имеет конечный ранг, если коразмерность векторного подпространства $P_\preceq\bigcap (-P_\preceq)$ является конечной, при этом коразмерность ${\rm codim}(P_\preceq\bigcap (-P_\preceq))$ будем называть \textit{рангом} $\preceq$.

%\begin{proposition}\label{pr4}
%Пусть $F$ --- грань выпуклого множества $Q \subset X$ и пусть согласованное полное отношение предпорядка $\preceq$, определенное на $X$, таково, что  $F = {\rm Min}(Q \mid \preceq) := \{ x \in Q \mid x \preceq y\,\,\text{для всех}\,\,y \in Q\}$. Тогда
%
%$($a$)$ грань $F$ является лексикографически выступающей в том и только том случае, когда на $X$ существует согласованное полное отношение предпорядка $\preceq$, ранг которого конечен, и такое, что $F = {\rm Min}(Q \mid \preceq) := \{ x \in Q \mid x \preceq y\,\,\text{для всех}\,\,y \in Q\}$. %коразмерность векторного подпространства $P_\prec \bigcap (-P_\prec)$ является конечной;
%
%$($b$)$ грань $F$ является выступающей в том и только том случае, когда на $X$ существует согласованное полное отношение предпорядка $\preceq$, ранг которого равен единице, и такое, что $F = {\rm Min}(Q \mid \preceq) := \{ x \in Q \mid x \preceq y\,\,\text{для всех}\,\,y \in Q\}$.%коразмерность векторного подпространства $P_\prec \bigcap (-P_\prec)$ равна единице, т.е. когда $P_\prec \bigcap (-P_\prec)$ является гиперплоскостью в $X$.
%\end{proposition}
%
%Доказательство Предложения \ref{pr4} будет дано ниже.

\section{The characterization of faces of convex sets by step-affine functions.}

First of all, let us recall the definition of step-affine functions and some their properties (more detailed information about step-affine functions is contained in \cite{GS2000,Gor20}).

Let $A(X,{\mathbb R})$ be the space of real-valued affine functions defined on a vector space $X$. Since ${\mathcal A}(X)=L(X)\times {\mathbb R},$ where $L(X)$ is the vector space of linear functions on $X$ and ${\mathbb{R}}$ is the set of real numbers, then the pair from $L(X)\times {\mathbb R}$ corresponding to an affine function $f\in {\mathcal A}(X)$ will be denoted by $(l_f,\alpha_f),$
with $f:x\rightarrow l_f(x)+\alpha_f.$

For any family $\mathcal{F} \subset \mathcal{A}(X)$, which is linearly ordered\footnote{A binary relation $\preccurlyeq$ defined on a set $Y$ is called a linear order relation if it is a total partial order relation.} by some relation $\preccurlyeq_{\mathcal{F}}$, and for any $f \in \mathcal{F}$ we define the set (affine manifold)
$$
E_f := \{x \in X \mid \bar{f}(x)=0\,\,\text{for all}\,\,\bar{f} \in \mathcal{F}\,\,\text{such that}\,\,\bar{f} \preccurlyeq_{\mathcal F} f,\, \bar{f} \ne f\},
$$

if $f$ is the smallest (first) element in $\mathcal{F}$, then we set $E_f = X$

We say that a family $\mathcal{F} \subset \mathcal{A}(X)$ forms \textit{a cortege of affine functions} on $X$ if

($i$) $\mathcal{F}$ is linearly ordered by some relation $\preccurlyeq_{\mathcal{F}}$;

($ii$) for every function $f \in \mathcal{F}$ the corresponding affine manifold $E_f$ is non-empty and $f(E_f)={\mathbb R};$

($iii$) for any $x \in X$ the subfamily ${\mathcal F}_x := \{f \in {\mathcal F} \mid f(x) \ne 0\}$ is either empty or has a least (with respect to $\preccurlyeq_{\mathcal{F}}$) element $f_x$.

If each function $f \in {\mathcal F}$ is in fact linear, then the cortege of affine functions ${\mathcal F}$ is called a \textit{cortege of linear functions}.

It follows from the property ($ii$) of the cortege ${\mathcal F}$ that $l_f \ne 0$ for any function $f \in {\mathcal{F}}$, i.e. the cortege of affine functions does not contain constant functions, and the cortege of linear functions does not contain the zero function. Moreover, the cortege of linear functions is \cite[Proposition~4.1]{GS2000} linearly independent in $L(X)$.

A real-valued function $u:X \to {\mathbb R}$ is called \textit{step-affine} (\textit{step-linear}) if there exists a cortege of affine (linear) functions ${\mathcal F}$ such that
\begin{equation}\label{e5.1}
u(x)= u_{\mathcal F}(x) :=\left\{
  \begin{array}{cr}
       0,&\text{when}\,\,{{\mathcal F}}_{x} = \varnothing,\\
     f_x(x),&\text{when}\,\,{{\mathcal F}}_{x} \ne \varnothing.\\
  \end{array}
     \right.
\end{equation}

If the cortege of affine functions ${\mathcal F}$ is finite, i.e. if ${\mathcal F} = \{f_1,f_2,\ldots,f_m\}$, and the linear ordering of ${\mathcal F}$ corresponds to the numbering, then the relation \eqref{e5.1} can be written as follows:
\begin{equation}\label{e5.1a}
u_{{\mathcal F}}(x)=\left\{
  \begin{array}{ll}
       f_1(x),&\text{when}\,\,f_1(x) \ne 0,\\
       f_2(x),&\text{when}\,\,f_1(x) = 0,f_2(x) \ne 0,\\
       ......... & ....................................................................... \\
       f_{m-1}(x),&\text{when}\,\,f_1(x) = \ldots = f_{m-2}(x)=0,f_{m-1}(x) \ne 0,\\
       f_{m}(x),&\text{when}\,\,f_1(x) = \ldots = f_{m-1}(x)=0.
       \end{array}
     \right.
\end{equation}

A step-affine function $u:X \to {\mathbb R}$ is called \textit{regular} if the set of points from $X$ at which it takes the zero value is non-empty, i.e. if $U_= := \{x \in X \mid u(x)=0\} \ne \varnothing$.

Since $U_= = \displaystyle\bigcap_{f \in {\mathcal F}}\{x \in X \mid f(x) = 0\}$, where ${\mathcal F}$ is a cortege of affine functions which defines (by \eqref{e5.1}) a step-affine function $u$, then the set $U_=$ is an affine manifold.

Every step-linear function is regular. Moreover, \cite[Proposition 4.4]{GS2000}, a step-affine function $u:X \to {\mathbb R}$ is regular if and only if there exist a step-linear function $w:X \to {\mathbb R}$ and a point $a \in X$ such that $u(x) = w(x -a)$ for all $x \in X$.

\begin{proposition}{\rm \cite[Theorem 4.1]{GS2000}}\label{pr5}
A proper subset $S$ of a vector space $X$ is a generalized semispace generated by an affine manifold $M \subset X$ if and only if there exists a regular step-affine function $u: X \to {\mathbb R}$ such that $S = \{x \in X \mid u(x) > 0\}$ and $M = \{x \in X \mid u(x) = 0\}$.
\end{proposition}

The statement of Proposition \ref{pr5} demonstrates the duality between regular step-affine functions and generalized semispaces generated by affine manifolds. It is easy to see that this duality is an extension of the classical duality between affine functions and algebraically open halfspaces.

\begin{theorem}\label{thC}
A subset $F \subset Q$ ia a proper face of a convex set $Q$ if and only if there exists a regular step-affine function $u:X \to {\mathbb R}$ defined on $X$ such that   $Q \subset \{x \in X \mid u(x) \geq 0\}$ and $F = \{x \in Q \mid u(x) = 0\}$.
\end{theorem}

\textit{Proof.}
\textit{Necessity.} Let $F$ be a proper face of a convex set $Q \subset X$. It follows from the necessary part of Theorem \ref{thA} that for $F$ there exist an affine manifold $M \subset X$ and a generalized semispace $S \subset X$, generated by $M$, such that $Q \subseteq M \bigcup S$ and $F = M \bigcap Q$. By Proposition \ref{pr5}, for these $S$ and $M$ there exists a regular step-affine function $u:X \to {\mathbb R}$ such that $S = \{x \in X \mid u(x) > 0\}$ and $M = \{x \in X \mid u(x) = 0\}$. It is easy to see that for such a function $u:X \to {\mathbb R}$, the inclusion $Q \subseteq M \bigcup S$ implies that $u(x) \ge 0$ for all $x \in Q$, and the equality $F = M \bigcap Q$ implies $F = \{x \in Q \mid u(x) = 0\}$.

\textit{Sufficiency.}
Let $u: X \to {\mathbb R}$ be a regular step-affine function such that $Q \subset \{x \in X \mid u(x) \geq 0\}$ and $F = \{x \in Q \mid u(x) = 0\}$ or, equivalently, such that $Q \subset S_u \bigcup M_u$ and $F = M_u \bigcap Q$, where $S_u:= \{x \in X \mid u(x) > 0\}$ and $M_u := \{x \in X \mid u(x) = 0\}$. Since, by Proposition \ref{pr5}, $M_u$ is a nonempty affine manifold in $X$, and $S_u$ is a generalized semispace generated by $M_u$, then, by the sufficient part of Theorem \ref{thA}, the subset $F$ is a proper face of the set $Q$. \hfill $\square$

It is clear that the characterizations of faces presented in Theorems \ref{thA}, \ref{thB}, \ref{thC} are equivalent to each other. We fix this fact explicitly in the following theorem.

\begin{theorem}\label{th4}
For any convex set $Q \subset X$ the following statements are equivalent:

$($a$)$ a subset $F \subset Q$ is a proper face of $Q$;

$($b$)$ there exist an affine manifold $M \subset X$ and a generalized semispace $S \subset X$, generated by $M$, such that
$Q \subset S \bigcup M$ and $F = Q \bigcap M$;

$($c$)$ there exists a compatible total preorder relation $\preceq$ defined on $X$ such that $F = {\rm Min}(Q \mid \preceq)$; %= {\rm Inf}(Q \mid \preceq)$;

$($d$)$ there exists a regular step-affine function $u: X \to X$ such that $u(x) \geq 0$ for all $x \in Q$ and $F = \{x \in Q \mid u(x)=0\}$.
\end{theorem}

\section{The characterization of lexicographically exposed faces}

We say that a step-affine function $u:X \to {\mathbb R}$ has \textit{finite rank} if the cortege of affine functions ${\mathcal F}$, which defines the function $u$ by \eqref{e5.1}, is finite, and the number of elements of the cortege ${\mathcal F}$ is called the \textit{rank} of the step-affine function $u$\footnote{The concept of rank for an arbitrary step-affine function is introduced in \cite{GS2000,Gor20}.}. A step-affine function whose rank is equal to one is in fact an affine function.

Note also, that a step-affine function $u_{\mathcal F}$ whose rank is finite is regular.

\begin{theorem}\label{th7}
Let $F$ be a proper face of a convex set $Q \subset X$.

$($a$)$ The face $F$ is lexicographically exposed if and only if on $X$ there exists a step-affine function $u:X \to {\mathbb R}$, whose rank is finite and which is   such that $Q \subset \{x \in X \mid u(x) \geq 0\}$ and $F = \{x \in Q \mid u(x) = 0\}$;
%коразмерность аффинного многообразия $M_u := \{x \in X \mid u(x) = 0\}$ конечна;
%ранг обобщенного семипространства $S$ конечен;

$($b$)$ the face $F$ is exposed if and only if on $X$ there exists a step-affine function $u:X \to {\mathbb R}$, whose rank is equal to one and which is such that $Q \subset \{x \in X \mid u(x) \geq 0\}$ and $F = \{x \in Q \mid u(x) = 0\}$.
%функция $u$ имеет конечный ранг равный единице, т.е. когда $u$ является аффинной функцией, отличной от константы.
%аффинного многообразия $M_u := \{x \in X \mid u(x) = 0\}$ равна единице, т.е. когда $M_u$ является гиперплоскостью в $X$.
\end{theorem}

\textit{Proof.}
($a$) \textit{Necessity.} Let $F$ be a lexicographically exposed face of a convex set $Q$. By Definition~1 there exists a finite, linearly independent, ordered family $\{l_1,l_2,\ldots,l_m\} \subset L(X)$ linear functions, defined on $X$, such that $F = \{x \in Q \mid l_i(x) = {\rm min}_{y \in F_{i-1}}l_i(y), i = 1,2,\ldots,m\},$ where $F_0 = Q, F_i = \{x \in F_{i-1} \mid l_i(x) = {\rm min}_{y \in F_{i-1}}l(y),\},  i = 1,2,\ldots,m$.
Define a family of affine functions ${\mathcal F} = \{f_1,f_2,\ldots,f_m\}$ by setting $f_i(x) := l_i(x) +\alpha_i,\,\,i=1,2,\ldots,m$, where $\alpha_i := - \min_{y \in F_{i-1}}l_i(y),\,\,i=1,2,\ldots,m$. Then $F = \{ x \in Q \mid f_i(x) =0, \,\,i=1,2,\ldots,m\}$. The family ${\mathcal F} = \{f_1,f_2,\ldots,f_m\}$, ordered according to the numbering of its elements, is a finite cortege of affine functions on $X$, and the step-affine function $u_{\mathcal F}$, which according to \eqref{e5.1a} is determined by the cortege ${\mathcal F} = \{f_1,f_2,\ldots,f_m\}$, has finite rank. It follows from \eqref{e5.1a} and the equalities $F_0 = Q, F_i = \{x \in F_{i-1} \mid f_i(x) = 0\}, i = 1,2,\ldots,m,$ that $Q \subset \{x \in X \mid u_{\mathcal F}(x) \geq 0\}$ and $F = \{x \in Q \mid u_{\mathcal F}(x) = 0\}$. Thus, the existence of the required step-affine function of finite rank that characterizes the lexicographically exposed face is proved.

\textit{Sufficiency.} Suppose that the step-affine function $u_{\mathcal F}$ defined by the finite cortege of affine functions ${\mathcal F} =\{f_1,f_2,\ldots,f_m\}$ is such that $Q \subset \{x \in X \mid u_{\mathcal F}(x) \geq 0\}$ and $F = \{x \in Q \mid u_{\mathcal F}(x) = 0\}$.
Since the step-affine function $u_{\mathcal F}$ is regular, then by Theorem \ref{thC} $F$ is a face of $Q$. We show that $F$ is, in fact, a lexicographically exposed face.
Recall that each affine function $f \in A(X)$ can be represented as $f(x) = l_f(x) + \alpha_f, x\in X$, where $l_f \in L(X), \alpha_f \in {\mathbb{R}}$.

%Покажем, что в этом случае множество $F = \{x \in Q \mid u_{\mathcal F}(x) = 0\}$ является лексикографически выступающей гранью выпуклого множества $Q$.
%
%Пусть $y,z \in Q$ и $\alpha \in (0,1)$ таковы, что $\alpha y + (1 - \alpha) z \in F$. Из \eqref{e5.1a} следует, что равенство $u_{\mathcal F}(x) = 0$ эквивалентно равенствам $f_1(x)=0,f_2(x)=0,\ldots,f_m(x)=0$. Так как функции $f_i, i = 1,2,\ldots,m$ аффинны, то для любого $i=1,2,\ldots,m$ имеем $0 = f_i(\alpha y + (1 - \alpha) z) = \alpha f_i(y) + (1 - \alpha) f_i(z)$. Поскольку $u_{\mathcal F}(x) \geq 0$ для любого $x \in Q$, то $f_i(y) \geq 0, f_i(z) \geq 0$ для всех $i = 1,2,\ldots, m$.  Следовательно, равенство $\alpha f_i(y) + (1 - \alpha) f_i(z) = 0$, где $\alpha \in (0,1)$, влечет $f_i(y) = 0$ и $f_i(z) = 0$. Значит, $y,z \in F$ Это доказывает, что $F$ является гранью $Q$.
Since $u_{\mathcal F}(x) \geq 0$ for all $x \in Q$, then $f_1(x) \geq 0$ for all $x \in Q$, which is equivalent to the inequality ${\rm min}_{y \in Q}l_{f_1}(y) + \alpha_{f_1} \geq 0$ and, consequently, the value ${\rm min}_{y \in Q}l_{f_1}(y)$ is finite. In turn, the equality $F = \{x \in Q \mid u_{\mathcal F}(x) = 0\}$ implies the inclusion $F \subseteq F_1 := \{ x \in Q \mid l_{f_1}(x) ={\rm min}_{y \in Q}l_{f_1}(y)\}$. If $F = F_1$, then the face $F$ is exposed, and therefore lexicographically exposed.

If $F_1 \ne F$, then from the condition $u_{\mathcal F}(x) \geq 0$ for all $x \in Q$ we conclude that $f_2(x) \geq 0$ for all $x \in F_1$, whence it follows that the quantity ${\rm min}_{y \in F_1}l_{f_2}(y)$ is finite. The equality $F = \{x \in Q \mid u_{\mathcal F}(x) = 0\}$ implies the inclusion $F \subseteq F_2 := \{ x \in F_1 \mid l_{f_2}(x) ={\rm min}_{y \in F_1}l_{f_2}(y)\}$. If $F = F_2$, then the face $F$
is lexicographically exposed.

If $F_2 \ne F$, then we continue constructing the sets $F_3 = \{ x \in F_2 \mid l_{f_3}(x) ={\rm min}_{y \in F_2}l_{f_3}(y)\}$, $F_4 = \{ x \in F_3 \mid l_{f_4}(x) ={\rm min}_{y \in F_3}l_{f_4}(y)\}$ and so on until at some step $k \leq m$ we arrive at the set $F_k = F$ and thereby establish that the face $F$ is lexicographically exposed. The fulfillment of the equality $F = F_k$ at some step $k \leq m$ is guaranteed by the equality $F = \{x \in Q \mid u_{\mathcal F}(x) = 0\}$.

($b$) Since a step-affine function whose rank is equal to one is in fact an affine function other than a constant, the statement ($b$) is equivalent to the definition of a exposed face.
\hfill $\square$

\smallskip

Since in finite-dimensional vector spaces any step-affine function has finite rank, then from Theorems \ref{thC} and \ref{th7} we arrive at the result previously established by J.-E. Martinez-Legaz \cite{ML}.

\begin{corollary}{\rm \cite{ML}.}
In finite-dimensional vector spaces, every proper face of any convex set is lexicographically exposed.
\end{corollary}

\section*{Conclusion}

Due to Proposition \ref{pr5} there is a dual correspondence between step-affine functions and generalized semispaces, using of which we can extend the notion of finite rank from step-affine functions to generalized semispaces. Moreover, since a compatible preorder is total if and only if $P_\prec := P_\preceq \setminus (-P_\preceq)$ is generalized semispace generated by the vector subspace $P_\preceq \bigcup (-P_\preceq)$, we can extend the notion of finite rank to compatible total preorders too. These extensions allow us to reformulate Proposition \ref{th7} in such way that we get the characterizations of lexicographically exposed faces both in the terms of generalized semispaces and in the terms of compatible total preorders having a finite rank. Every of these characterizations implies that in finite-dimensional spaces each face of a convex set is lexicographically exposed.

Thus, the characterizations of faces of convex sets presented in Theorems \ref{thA}, \ref{thB} and \ref{thC} extend to infinite-dimensional vector spaces the lexicographical characterization of faces proved by Martinez-Legaz J.-E. in \cite{ML} for finite-dimensional vector spaces.

\section*{Acknowledgements}
The research was supported by the State Program for Fundamental Research of Republic of Belarus.
%\section*{Compliance with ethical standards}

\section*{Conflict of Interest}
The author declares that he have no conflict of interest.


\begin{thebibliography}{99}

\bibitem{Rock}
Rockafellar, R.T.: Convex Analysis, Princeton University Press, Princeton, New Jersey (NJ)  (1973).

\bibitem{Boris}
Mordukhovich, B.S., Nam, N.M.: Convex analysis and beyond. Vol.1 Basic theory. Springer (2022)  https://doi.org/10.1007/978-3-030-94785-9; ISBN 9783030947842

\bibitem{Bronsted}
Bronsted, A.: An Introduction to Convex Polytopes. Springer-Verlag (1983) https://doi.org/10.1007/978-1-4612-1148-8; ISBN 978-0-387-90722-2

\bibitem{Millan}
Diaz Millan, R., Roshchina V.: The Intrinsic Core and Minimal Faces of Convex Sets in General Vector Spaces. Set-Valued and Variational Analysis, Article 14 (2023) https://doi.org/10.1007/s11228-023-00671-6 %27~p.


\bibitem{ML}
Martinez-Legaz, J.-E.: Lexicographical characterization of the faces of convex sets. Acta Mathematica Vietnamica 22(1), 207--211 (1997)


\bibitem{Las}
Lassak, M.: Convex half-spaces. Fund. Math. 120(2), 7--13 (1984) https://doi.org/10.4064/fm-120-1-7-13 %21

\bibitem{GS98}
Gorokhovik, V.V., Semenkova, E.A.: Classification of semispaces according to their types in infinite-dimensional vector spaces. Mathematical Notes 64(2), 164 -- 169  (1998) https://doi.org/10.1007/BF02310300

\bibitem{GS2000}
Gorokhovik, V.V., Shinkevich, E.A.: Geometric structure and classification of infinite-dimensional halfspaces. In: Przeworska--Rolewich, D. (ed). Algebraic Analysis and Related Topics. Banach Center Publications 53, 121--138 (2000)  https://www.researchgate.net/publication/388932409


\bibitem{Gor20}
Gorokhovik, V.V.: Step-Affine Functions, Halfspaces, and Separation of Convex Sets with Applications to Convex Optimization Problems.  Proceedings of
the Steklov Institute of Mathematics (Suppl.) 313(Suppl.~1), S83--S99 (2021) https://doi.org/10.1134/S008154382103010X %8

\bibitem{Hille}
Hille, E., Phillips, R.S.: Functional analysis and Semi-Groups. Pronidence (RI): American Math. Society (1957) %819~p.

\bibitem{Jech}
Jech, T.J.: Set theory: The third millenium edition (3rd ed.). Springer-Verlag (2006) %787~p.(ISBN 3540440852)

\bibitem{Peressini}
Peressini, A.L.: Ordered topological vector spaces. Harper and Row, (1967)

\end{thebibliography}
\end{document}